# Fibre dimension of the Nash transformation

Achim Hennings[♦]

**Abstract**

It is proved that the special fibre of the Nash transform of an irreducible isolated singularity has maximal possible dimension.

In the paper [2] by A. Simis, K. Smith and B. Ulrich it is asked (in slightly different terms) as a question, whether at the Nash transformation of an isolated singularity the special fibre has always the maximal dimension. In the particular case of a Cohen-Macaulay singularity a positive answer is given there, using a result of E. Kunz and R. Waldi [1]. Here, we propose a simple proof for the general case.

Let $(X, 0) \subseteq (\mathbb{C}^N, 0)$ be a reduced isolated singularity of dimension $d \geq 2$. We suppose $X$ to be irreducible and non-smooth. The smooth locus will be denoted by $X^*$.

The Nash transform $\nu: \hat{X} \to X$ can be described geometrically as the closure of the graph of the Gauss map $X^* \to G$, $x \mapsto T_{X,x}$, which assigns to every smooth point its tangent space considered as a point of the Grassmann manifold of $d$−dimensional subspaces of $\mathbb{C}^N$. In this way $\hat{X} \subseteq X \times G$ is an analytic subset, and $\nu$, which is induced from the first projection, is a proper modification. The special fibre $\nu^{-1}(0) \subseteq G$ is just the set of all limits of tangent spaces at smooth points. The universal bundle on $G$ provides $\hat{X}$ with a (dual) locally free sheaf $\widetilde{\Omega}$ and a surjection $\nu^*\Omega_X^1 \to \widetilde{\Omega}$, which is isomorphic off $\nu^{-1}(0)$. (Cf. e.g. [3])

Let now $n: Y \to \hat{X}$ be the normalisation and $\pi: Y \to X$ the composition with the Nash transformation. The sheaf $\widetilde{\omega} \coloneqq n^* \Lambda^d \widetilde{\Omega}$ is invertible and outside $\pi^{-1}(0)$ it coincides with the canonical sheaf $\omega_Y$. We note that for a normal complex space the canonical sheaf is the sheaf of all forms of top degree which need only be defined and regular at smooth points.

**Theorem:** The fibre $\pi^{-1}(0)$ has pure dimension $d - 1$. In particular the same holds for the fibre $\nu^{-1}(0)$ of $\hat{X}$.

Proof: By its construction as a closure of a manifold of dimension $d$ the Nash transform is pure dimensional. This property is inherited by the normalisation $Y$. Therefore it is sufficient to disprove, that there is a point $P \in \pi^{-1}(0)$ at which $\pi^{-1}(0)$ has codimension $\geq 2$. Under this assumption however, taking into account that $\widetilde{\omega}$ is invertible, we must have equality $\widetilde{\omega}_P = \omega_{Y,P}$. This means that we have a surjective map $(\pi^*\Omega_X^d)_P \to \omega_{Y,P}$. Our aim is to show that this is impossible.

Case 1: $P$ is a smooth point in $Y$. The map $\pi: Y \to X \to \mathbb{C}^N$ has in local coordinates $(y_1, \ldots, y_d)$ around $P$ the form $\pi: (\mathbb{C}^d, 0) \to (\mathbb{C}^N, 0)$, and the image of the map

---

[♦] Universität Siegen, Fakultät IV, Hölderlinstraße 3, D-57068 Siegen

$$(\pi^*\Omega_X^d)_P \to \omega_{Y,P} \cong \mathcal{O}_{Y,P}$$

is the ideal of maximal minors $I_d(\partial\pi/\partial y)$ of the Jacobian. From this fact we see that the map can only be surjective when $\pi$ is regular at $P$. But then $X$ would have a smooth component.

Case 2: $P$ is a singular point in $Y$. We consider the decomposition:

$$(\pi^*\Omega_X^d)_P \to \Omega_{Y,P}^d \to \omega_{Y,P}$$

The map at right is of course also surjective. To get a contradiction, let $y_1, \ldots, y_m$ be a general linear coordinate system on some embedding space $(\mathbb{C}^m, 0) \supseteq (Y, P)$. As $\omega_{Y,P} \cong \mathcal{O}_{Y,P}$ and $\Omega_{Y,P}^d$ is generated by all $dy_{i_1} \wedge \ldots \wedge dy_{i_d}$ $(1 \leq i_1 < \cdots < i_d \leq m)$, we can suppose that $\omega_{Y,P} = \mathcal{O}_{Y,P}(dy_1 \wedge \ldots \wedge dy_d)$. By the general choice of coordinates, the mapping

$$(y_1, \ldots, y_d) : (Y, P) \to (\mathbb{C}^d, 0)$$

is finite. However, this map has an empty ramification divisor. $Y$ being normal, this is impossible.

**Remark:** The above proof is also valid for non-isolated singularities. In this case one obtains a lower estimate of the fibre dimension by $d-1$ minus the dimension of the singular locus. This can also be obtained by taking hyperplane sections.